[1994/12/01]% LaTeX date must December 1994 or later
\documentclass[11pt,reqno]{article}

\title{The critical Branching Markov Chain is transient}
\author{Nina Gantert and  Sebastian M\"uller}
\date{25.10.05}
\usepackage{amsmath,amsthm,amsfonts}
\usepackage{t1enc}
\chardef\bslash=`\\ % p. 424, TeXbook
\newtheorem{thm}{Theorem}[section]

\newtheorem{lem}[thm]{Lemma}

\theoremstyle{definition}
\newtheorem{defn}{Definition}[section]

\theoremstyle{remark}
\newtheorem{rem}{Remark}[section]

\newcommand{\eps}{\varepsilon}

\newcommand{\R}{\mathbb{R}}
\newcommand{\Z}{\mathbb{Z}}
\newcommand{\N}{\mathbb{N}}
\renewcommand{\P}{\mathbb{P}}

\newcommand{\eval}[2][\right]{\relax
  \ifx#1\right\relax \left.\fi#2#1\rvert}

\begin{document}
\maketitle \abstract We investigate recurrence and transience of Branching Markov Chains
(BMC) in discrete time. Branching Markov Chains are clouds of particles which move (according to an irreducible underlying Markov Chain) and produce offspring
 independently. The offspring distribution can depend on the location of the particle. If the offspring distribution is constant for all locations, these are Tree-Indexed Markov chains in the sense of \cite{benjamini94}.
Starting with one particle at location $x$, we denote by
$\alpha(x)$ the probability that $x$ is visited infinitely often
by the cloud. Due to the irreducibility  of the underlying Markov
Chain, there are three regimes: either $\alpha(x) = 0$ for all $x$
(transient regime), or $0 < \alpha(x) < 1$ for all $x$ (weakly
recurrent regime) or $\alpha(x) = 1$ for all $x$ (strongly
recurrent regime). We give classification results, including a
sufficient condition for transience in the general case.
 If the mean of the offspring distribution is constant,
we give a criterion for transience involving the spectral radius of the underlying Markov Chain and the
 mean of the offspring distribution.
In particular, the critical BMC  is transient.
Examples for the classification
are provided.
\newline {\scshape Keywords:} Branching Markov Chains, recurrence and transience,
Lyapunov function,  spectral radius
\newline {\scshape AMS 2000 Mathematics Subject Classification:}
60J10, 60J80
\renewcommand{\sectionmark}[1]{}

\section{Introduction}
A Branching Markov Chain (BMC)  is a system of particles, which
move independently according to the transition probabilities of an
underlying Markov chain. We take a countable state space $X$
and an irreducible stochastic transition kernel $P$ for the
underlying Markov chain $(X,P).$  The BMC starts with one particle
in an arbitrary starting position $x_s \in X$ at time $0$.
Particles move independently according to $P.$  At each position
$x \in X$, they independently produce offspring according to some
probability distribution $\mu(x)$ on $\{1,2,3, \ldots \}$ (which
can depend on the position $x$ of the particle) and die. We  assume that
there is always at least one offspring particle, so that the
number of particles is always increasing in time. Similar models
have been
studied in \cite{menshikov97}.\\
The transition probabilities of the Markov chain and the offspring
distribution can be given as a (typical) realization of a random
environment. The behavior of the resulting \lq\lq Branching Random
Walk in Random Environment\rq\rq\ has been classified in
\cite{comets98}, \cite{machado00} and \cite{machado03} for the
case where the underlying Markov chain is a Random Walk in Random
Environment on $\Z^+$
or on a tree. A similar,
but more general model, where movement and offspring production are not independent anymore,
 is considered in \cite{comets05}. \\

Let $\alpha(x)$ be the probability that, starting the BMC from
$x_s =x$, the location $x$ is visited by infinitely many
particles. Using the irreducibility of the underlying Markov Chain, we obtain, similar to
Lemma 3.1 in \cite{benjamini94}, the following classification:

\begin{lem}\label{class}
There are three possible regimes:
\begin{equation}\label{trans}
\alpha(x) = 0\quad \forall x\in X
\end{equation}
(transient regime)
\begin{equation}\label{between}
0 < \alpha(x) < 1  \quad \forall x\in X
\end{equation}
(weakly recurrent regime)
\begin{equation}\label{stronglyrec}
\alpha(x) = 1\quad \forall x\in X
\end{equation}
(strongly recurrent regime).
\end{lem}
We write $\alpha=0$ $(>0,=1)$ if $\alpha(x) = 0~ (>0,=1)~ \forall
x\in X.$ We say that a BMC is {\it recurrent} if it is not
transient, i.e. if (\ref{between}) or (\ref{stronglyrec}) are satisfied. Note that in the weakly recurrent regime, the values
of $\alpha(x)$ do in general not coincide.

We first give a sufficient condition, Theorem \ref{thm:1}, for
transience where the Markov chain can be any irreducible Markov
chain and the branching distributions can be arbitrary. Under the
assumption of constant mean offspring we obtain in Theorem
\ref{thm:0} a classification in transience and recurrence for all
irreducible Markov chains. In particular, we show that in the
critical case the BMC is transient. It is left to forthcoming work
to study the subdivision of the recurrent phase. Under
homogeneity conditions, i.e. quasi-transitivity,  on the BMC we
show that the strongly recurrent regime coincides with the recurrent regime, i.e. (\ref{between}) does not occur, see Theorem
\ref{thm:2}.

\section{Preliminaries}

We give the definition of the spectral radius of an irreducible
Markov chain $(X,P)$ and quote a result which characterizes the
spectral radius in terms of $t-$superharmonic functions. For
further details see e.g. \cite{woess}.
\begin{defn}
Let $(X,P) $ be an irreducible Markov chain with countable state
space $X$ and transition operator $P=\left( p(x,y)\right)_{x,y\in
X}.$ The spectral  radius of $(X,P)$ is defined as
\begin{equation}
\rho(P):=\limsup_{n\rightarrow\infty} p^{(n)}(x,y)^{1/n} \in
(0,1],
\end{equation}
where $p^{(n)}(x,y)$ is the probability to get from $x$ to $y$ in
exactly $n$ steps. $P$ is interpreted as a (countable) stochastic
matrix, so that $p^{(n)}(x,y)$ is the $(x,y)-$entry of the matrix
power $P^n.$ We set $P^0=I,$ the identity matrix over $X.$
\end{defn}
The transition operator $P$ acts on functions $f:~X\rightarrow\R$
by
\begin{equation}
Pf(x):= \sum_y p(x,y) f(y).
\end{equation}
\begin{defn}
The Green function of $(X,P)$ is the power series
$$G(x,y|z)=\sum_{n=0}^\infty p^{(n)} (x,y) z^n,~ x,y\in X,~
z\in\mathbb{C}.$$
\begin{rem}\label{rem:2a}
For all $x,y\in X$ the power series $G(x,y|z)$ has the same radius
of convergence $1/\rho(P).$
\end{rem}
\end{defn}

\begin{defn} Fix $t>0.$
A $t-$superharmonic function  is a function $f:~X\rightarrow\R$
satisfying
$$Pf\leq t f.$$
We write $S(P,t)$ for the collection of all $t-$superharmonic
functions and $S^+(P,t)$ for the positive cone of $S(P,t)$, i.e.
$S^+(P,t) = \{f \in S(P,t): f \geq 0\}$.
\end{defn}

 A base of the
cone $S^+(P,t)$ can be defined with the help of a reference point
$x_0\in X$ by
$$B(P,t):=\{f\in S^+(P,t):~ f(x_0)=1\}.$$

\begin{lem}\label{lem:1a}
$B(P,t)$ is compact in the topology of
pointwise convergence.

\end{lem}
\begin{proof}
The closedness of $B(P,t)$ follows from Fatou's lemma. Let $x\in
X,$ then irreducibility implies the existence of  $n_x$ such that
$p^{(n_x)}(x_0,x)>0.$ If $f\in B(P,t)$ then
$$p^{(n)}(x_0,x)f(x)\leq P^n f(x_0) \leq t^nf(x_0)=t^n.$$ Hence
$$f(x)\leq \frac{t^{n_x}}{p^{(n_x)}(x_0,x)} \quad \forall f\in
B(P,t),$$ and the desired compactness follows.
\end{proof}

\begin{lem}\label{lem:1}
$$\rho(P)=\min\{t>0:~\exists\, f(\cdot)>0\mbox{ such that }Pf\leq t f\}$$
\end{lem}
\begin{proof}
If there exists a function $f\neq 0$ in $S^+(P,t),$ then
$p^{(n)}(x,x)f(x) \leq P^n f(x)\leq t^n f(x).$ Hence
$\rho(P)=\limsup_n p^{(n)}(x,x)^{1/n}\leq t.$ Conversely, for
$t>\rho(P)$ the function $f(x)=G(x,x_0|1/t)$ is by Remark
\ref{rem:2a} well-defined. It is clear that $f(\cdot)$ is non-zero
and in $S^+(P,t).$ Hence, $B(P,t)\neq\emptyset.$ We have
$B(P,t_1)\subseteq B(P,t_2)$ for $t_1<t_2.$ By compactness of the
sets $B(P,t),$ it follows that $B(P,\rho(P))=\bigcap_{t>\rho(P)}
B(P,t)\neq\emptyset.$
\end{proof}

\subsection{Branching Markov Chains}\label{sect:BMC}
We consider an irreducible Markov chain $(X,P)$ in discrete time.
For all $x\in X$ let
$$\mu_1(x),\mu_2(x),\ldots$$ be a sequence of non-negative numbers
satisfying
$$\sum_{k=1}^\infty \mu_k(x)=1 \mbox{ and } m(x):=\sum_{k=1}^\infty k \mu_k(x)<\infty.$$
  We define the Branching Markov Chain (BMC) on $(X,P)$ following \cite{menshikov97}.
  At time $0$ we start with one particle in an arbitrary starting position
$x_s\in X.$ When a particle is in $x$, it generates $k$ offspring
particles at $x$ with probability $\mu_k(x)$ ($k=1,2, \ldots $) and dies. The $k$
offspring particles then move independently according to the
Markov chain $(X,P)$ and generate their offspring as well. At any
time, all particles move and branch  independently of the other
particles and the previous history of the process. The resulting
BMC is a Markov chain with countable state space $X'$, namely the
space of all particle configurations
$$\omega(n)=\{x_1(n),x_2(n),\ldots,x_{\eta(n)}(n)\},$$
where $x_i(n)\in X$ is the position of the $i$th particle at time $n$ and
$\eta(n)$ is the total number of particles at time $n$. Since
there is always at least one offspring particle, the number of particles is
always increasing in time. In most cases under consideration the
number of particles $\eta(n)$ tends to infinity as
$n\rightarrow\infty$ almost surely. Therefore, it is not
interesting to ask if a BMC is recurrent as a Markov chain
on $X'$: $\eta(n)\rightarrow\infty$ implies its transience. It is
more reasonable to define transience and recurrence as in Lemma
\ref{class}. With the notations above we can write $\alpha$
as
$$ \alpha(x)=\P_x\left( \sum_{n=1}^\infty
\sum_{i=1}^{\eta(n)}\mathbf{1}_{\{x_i(n)=x\}}=\infty\right),$$ where
$\P_x(\cdot)=\P(\cdot|x_s=x)$ and $x\in X.$
 Note that a BMC in our setting is strongly recurrent ($\alpha=1$) if
every state $x\in X$ is visited with probability 1. In analogy to
\cite{menshikov97}, we introduce the following modified version of
the BMC. We fix an arbitrary position $x_0\in X,$ which we denote
the origin of $X.$ After the first time step we conceive the
origin as an {\it absorbing} point: if a particle reaches the
origin it stays there forever and stops producing offspring. We
denote  this new process with BMC*. The process BMC* is analogous
to the original process BMC except that $p_{x_0,x_0}=1,~
p_{x_0,x}=0~\forall x\neq x_0$ and $\mu_1(x_0)=1$ from the second
time step on. Let $\eta_0(n,x_s)$ be the number of particles at
position $x_0$ at time $n$, given that the BMC* started in $x_s\in
X.$ We define the random variable $\nu(x_s)$ as
$$\nu(x_s)=\lim_{n\rightarrow\infty} \eta_0(n,x_s).$$
The random variable $\nu$ takes values in
$\{0,1,2,\ldots\}\cup\{\infty\}.$

\section{Results} We present a sufficient condition for
transience of a Branching Markov Chain (BMC), which is inspired by
the Lyapunov methods developed in \cite{comets98} and
 \cite{menshikov97}.

\begin{thm}\label{thm:1}
A BMC with irreducible underlying Markov chain $(X,P)$ and
$m(y)>1$ for some $y\in X$ is transient if there exists a strictly
positive function $f(\cdot)$ such that
\begin{equation}\label{eq:1.1}
P f(x)\leq \frac{f(x)}{m(x)}\quad \forall x\in X.
\end{equation}
\end{thm}

\begin{proof}
We show that the total number of particles returning to a starting
point $x_s=x_0\neq y$ is finite. The total number of particles in
$x_0$ can be interpreted as the total number of progeny in a
branching process $(Z_n)_{n\geq 0}$. We show that this process
dies out with probability one. The branching process $(Z_n)_{n\geq
0}$ is defined as follows: Note that each particle has a unique
ancestry line which leads back to the starting particle at time
$0$ at $x_0$. Let $Z_0=1$ and let $Z_1$ be the number of particles
being the first particle in their ancestry line (after the
starting particle) to visit $x_0$. Inductively we define $Z_n$ as
the number of particles being the $n$th particle in their ancestry
line to visit $x_0$. This defines a Galton-Watson process with
offspring distribution $Z{\buildrel d \over =} Z_1.$ We have that
$$\sum\limits_{n=1}^\infty Z_n =\sum\limits_{n=1}^\infty \sum_{i=1}^{\eta(n)}\mathbf{1}_{\{x_i(n)=x_0\}}$$
 and $Z{\buildrel d \over =} \nu(x_0).$
In order to show that $(Z_n)$ dies out almost sure it suffices to
show that $E\nu(x_0)\leq 1$ and $\P_{x_0}(\nu(x_0)<1)>0.$ Given the first
statement the latter is true since $m(y)>1$ and hence $\P_{x_0}(\nu(x_0)>
1)>0.$ It remains to show the first statement:
 Consider the corresponding BMC* and define
$$Q(n):=\sum_{i=1}^{\eta(n)} f(x_i(n)),$$
where $x_i(n)$ is the position of the $i$th particle at time $n.$
One can show that $Q(n)$ is a supermartingale. We refer the reader
for the technical details to the proof of Theorem 3.2 in
\cite{menshikov97}.

As $Q(n)$ is a positive supermartingale it converges almost surely
to a random variable $Q_\infty.$ Fatou's Lemma implies
$$EQ_\infty\leq \lim_{n\rightarrow\infty} EQ(n) \leq EQ(0).$$
For a BMC* started in a position $x_s\in X$ we also have  that
$$Q(n)\geq \eta_0(n,x_s) f(x_s)$$ and hence that
$$\nu(x_s)\leq \frac{Q_\infty}{f(x_s)}.$$ We obtain by taking
expectations and starting the BMC* in $x_s=x_0$
\begin{equation}\label{eq:1.11}
E\nu(x_0)\leq
\frac{EQ_\infty}{f(x_0)}\leq\frac{EQ(0)}{f(x_0)}=\frac{f(x_0)}{f(x_0)}=1.\end{equation}
\end{proof}
\begin{rem}
In contrast to Theorem 2.2. in \cite{comets98} and Corollary 3.1
in \cite{menshikov97} we demand that the condition (\ref{eq:1.1})
holds for all $x\in X.$ Note that in \cite{comets98} the BMC* is
defined in a slightly different way: the origin $x_0$ is always
absorbing.
\end{rem}

\begin{rem}
The converse of Theorem \ref{thm:1} does not hold in
general, for a counterexample see Section 5 in \cite{comets98}.
\end{rem}

\subsection{BMC with constant mean offspring}
 We assume that the mean number of offspring is constant, i.e.
 $m(x)=m>1$ for all $x\in X.$ Note that we do
not assume $(\mu_k(x))_k=(\mu_k(y))_k$ for $x,y\in X$, and the BMC
therefore needs not to be a Tree-Indexed Markov Chain as in
\cite{benjamini94}.

Under these assumptions, we have the following.

\begin{thm}\label{thm:0}
For a BMC with irreducible underlying Markov chain $(X,P)$ and
constant mean offspring $m>1$, it holds that the BMC is transient
if $m\leq 1/\rho(P)$ and recurrent if $m>1/\rho(P).$
\end{thm}
\begin{rem}
If $m=\infty$ then the BMC is  recurrent, since one can compare
the process with a suitable BMC with $\widetilde{m}>1/\rho(P).$
\end{rem}
\begin{proof}
The first part follows from Lemma \ref{lem:1} and Theorem
\ref{thm:1}. To show the recurrence we use ideas developed in
\cite{benjamini94} and \cite{comets98}: In order to show the
recurrence we compare the original BMC by some new process with
fewer particles and show that this process is recurrent. We start
the BMC in $x_0\in X$. We know from the hypothesis and the
definition of $\rho(P),$ that there exists a $k=k(x_0)$ such that
$$p^{(k)}(x_0,x_0)> m^{-k}.$$ We construct a new process $\xi(\cdot)$ by observing the BMC
only  at times $k, 2k, 3k, \ldots $ and by killing all the
particles not being in position $x_0.$ Let $\xi(n)$ be the number
of particles of the new process in $x_0$ at time $nk.$  The
process $\xi(\cdot)$ is a Galton-Watson process with mean
$p^{(k)}(x_0,x_0)\cdot m^{k} > 1,$ thus survives with positive
probability and hence the origin is hit infinitely often with
positive probability.
\end{proof}

\begin{rem}
Theorem \ref{thm:0} implies in particular that Markov chains
indexed by Galton-Watson trees are transient in the critical case
$m=1/\rho(P),$ since if $(\mu_k(x))_k=(\mu_k(y))_k$ for all
$x,y\in X$ the BMC is a  Markov chain indexed by  a Galton-Watson
tree, compare to \cite{benjamini94}.
\end{rem}

\subsection{Quasi-transitive BMC}
Let $X$ be a locally finite, connected graph and $Aut(X)$ be the
group of automorphisms of $X.$ Let $P$ be the transition matrix of
an irreducible  random walk on $X$ and $Aut(X,P)$ be the group of
all $\gamma\in Aut(X)$ which satisfy $p(\gamma x,\gamma y)=p(x,y)$
for all $x,y\in X.$ We say the Markov chain $(X,P)$ is transitive,
if the group $Aut(X,P)$ acts transitively on $X$ and
quasi-transitive if $Aut(X,P)$ acts with finitely many orbits on
$X,$ that is that each vertex of $X$ belongs to one of finitely
many orbits.

We say a BMC is quasi-transitive if the group $Aut(X,P,\mu)$ of
all $\gamma\in Aut(X,P)$ which satisfy $\mu_k(x)=\mu_k(\gamma
x)~\forall k\geq 1$ for all $x\in X$ acts with finitely many
orbits on $X.$ Using induction on $n$, one can show the following.

\begin{lem}\label{lem:quasibmc}
For a quasi-transitive BMC it holds that for all $x,y \in X$ and
all $\gamma\in Aut(X,P,\mu)$

\begin{equation}\label{eq:lem:quasibmc}
\P_x\left(\sum_{i=1}^{\eta(n)} \mathbf{1}\{x_i(n)=y\}=k\right)=
\P_{\gamma x}\left(\sum_{i=1}^{\eta(n)} \mathbf{1}\{x_i(n)=\gamma
y\}=k\right)\quad \forall n\in \N.
\end{equation}
\end{lem}

For quasi-transitive BMC we have a $0-1-$ law for the return
probability. In other words, $\alpha\in\{0,1\}$ in this case.

\begin{thm}\label{thm:2}
For a quasi-transitive BMC with underlying Markov chain $(X,P)$
and branching distribution $(\mu_k(x))_{k\geq 1}$ with constant
mean offspring $m(x)= m  >1, \forall x$, it holds that
\begin{itemize}
    \item the BMC is transient $(\alpha=0)$ if $m\leq 1/\rho(P)$.
    \item the BMC is strongly recurrent $(\alpha=1)$ if $m>1/\rho(P)$.
\end{itemize}
\end{thm}

\begin{proof}
The statement for the case $m\leq 1/\rho(P)$ follows from Theorem \ref{thm:0}. Recurrence in the case $m>1/\rho(P)$ also follows from Theorem \ref{thm:0}. In order to show
the strong recurrence ($\alpha=1$) in the case $m>1/\rho(P)$, we have to refine the arguments
from the proof of Theorem \ref{thm:0}. Constructing infinitely
many supercritical Galton-Watson processes whose extinction
probabilities are bounded away from $1$, we show that at least one
location is hit infinitely often.  We start the  BMC in
$x_{s_1}\in X$. We know from the hypothesis and the definition of
$\rho(P),$ that there exists a $k_1=k_1(x_{s_1})$ such that
$$p^{(k_1)}(x_{s_1},x_{s_1})> m^{-k_1}.$$
We construct a new process $\xi_1(\cdot)$ by observing the BMC
only at times $k_1, 2k_1, 3k_1, \ldots $ and by killing all the particles not
being in position $x_{s_1}$. Then, $\xi_1(n)$ is the number of particles of the new process in $x_{s_1}$ at time $nk_1$. In this way, we obtain a Galton-Watson process $\xi_1(\cdot)$. The number of
particles in $x_{s_1}$ at time $nk_1$ of the original BMC is at least $\xi_1(n)$. The process $\xi_1(\cdot)$
is a Galton-Watson process with mean
$p^{(k_1)}(x_{s_1},x_{s_1})\cdot m^{k_1} > 1.$ Hence
$\xi_1(\cdot)$ dies out with a probability $q_1=q_1(x_{s_1})<1.$
If this first process dies out, we start a second process
$\xi_2(\cdot)$, defined in the same way with a starting position $x_{s_2}$ ($x_{s_2}$ can be
any location which is occupied by a particle at the time where the first process dies out) and
$k_2=k_2(x_{s_2})$ such that
$$p^{(k_2)}(x_{s_2},x_{s_2})> m^{-k_2}.$$ This process dies out with
probability $q_2=q_2(x_{s_2}).$ If the second process dies out we
construct a third one, and so on. We obtain a sequence of processes $\xi_i(\cdot)$
with extinction probabilities $q_i.$ It suffices now to show that
the $q_i$ are bounded away from $1$: the probability
that all the processes die out is then $\prod_i q_i=0.$ Due to Lemma
\ref{lem:quasibmc} we have that for two starting positions $x$ and $y$
of the same orbit

$$\P_x\left(\sum_{i=1}^{\eta(n)} \mathbf{1}\{x_i(n)=x\}=k\right)=
\P_{y}\left(\sum_{i=1}^{\eta(n)}
\mathbf{1}\{x_i(n)=y\}=k\right)\quad \forall n\in \N.
$$
Hence two processes started in $x$ and $y$ have the same distributions
and hence the same extinctions probabilities. Since there are only
finitely many orbits, there are only finitely many different extinction
probabilities $q_i$.
\end{proof}

\begin{rem} Instead of considering quasi-transitive Markov Chains, we could also assume that
$(p^{(l)}(x,x))^{1/l}$ converges uniformly in $x$, i.e.
$\forall~\eps>0~\exists~l:~
(p^{(l)}(x,x))^{1/l}>\rho(P)-\eps~\forall x\in X,$
 and that there is a $k\in\N$ such that
$\inf_x \sum_{i=1}^k i \mu_i(x) \geq 1/\rho(P).$ Observing in the
same way as in the proof of Theorem \ref{thm:2} the BMC with
branching distributions $\tilde\mu_0(x)=\sum_{i=k+1}^\infty
\mu_i(x)$ and $\tilde\mu_i(x)=\mu_i(x)$ for $i=1,\ldots,k$ and
$x\in X,$ we obtain supercritical Galton-Watson processes $\xi_i$
with bounded variances and means bounded away from $1$, since $l$
and $k$ do not depend on $x_{s_i}.$ Hence the extinction
probabilities $q_i$ are bounded away from $1.$
\end{rem}

\section{Examples}
\begin{enumerate}
    \item A BMC with transient underlying Markov chain $(X,P)$ is
    transient if $$\sup_{x\in X} m(x)\leq 1/\rho(P).$$
    \item  A branching symmetric random walk on $\Z^d,$ $d\in\N,$
    is strongly recurrent for all branching distributions with constant mean offspring $m >1$.
    \item Consider a random walk on $\Z$ with drift: Let
    $X=\Z$, $p\in(0,1)$ and $P$ given by
    $$p(x,x+1)=p=1-p(x,x-1).$$ Take branching distributions with constant mean offspring $m$.
    The spectral radius is $\rho(P)=2\cdot \sqrt{p(1-p)}.$ Hence, the
    corresponding BMC is transient if  $$m\leq
    \frac1{2 \cdot\sqrt{p(1-p)}}$$  and strongly recurrent if
    $$m>
    \frac1{2 \cdot\sqrt{p(1-p)}}\, .$$
    (This reproduces a result of
    \cite{comets98} in section $4,$ noted that there is a
    calculation error in the formula after Theorem 4.3 of \cite{comets98} so that
    the "$<$" should become a "    $\leq$".)
    \item More generally, take $X=\Z^d$ and $e_i\in\Z^d$ with $(e_{i})_j=\delta_{ij}$
    for $i,j\in\{1,\ldots,d\},~ d\geq 1.$ Let $P$ be defined by
    $$ p(x,x+e_i)=p_i^+,~ p(x,x-e_i)=p_i^- ~\mbox{ such that } $$ $$\sum_{i=1}^d p_i^+
    +\sum_{i=1}^d p_i^-=1,\quad\forall x\in\Z^d$$ and such that $P$ is
    irreducible. Take branching distributions with constant mean offspring $m$.
    The spectral radius can be calculated with the help of the
    Perron-Frobenius Theorem (see for example \cite{woess}):
    $$\rho(P)=2\sum_{i=1}^d \sqrt{p_i^+ p_i^-}.$$
    The corresponding BMC is strongly recurrent  if
    $$ m >\frac1{2\sum_{i=1}^d \sqrt{p_i^+ p_i^-}}\, .$$
    Otherwise it is transient.
\end{enumerate}

\bigskip
 \noindent
\begin{tabular}{l}
 Nina Gantert\\
 Institut f\"ur Mathematische Statistik\\
 Universit\"at Münster\\
 Einsteinstr. 62  \\
 D-48149 Münster\\
 Germany\\
{\tt gantert@math.uni-muenster.de}\\
\\ \\
Sebastian M\"uller \\
 Institut f\"ur Mathematische Statistik\\
 Universit\"at Münster\\
 Einsteinstr. 62  \\
 D-48149 Münster\\
 Germany\\
{\tt Sebastian.Mueller@math.uni-muenster.de}\\
\end{tabular}

\end{document}